\documentclass[]{interact}

\usepackage{epstopdf}
\usepackage[caption=false]{subfig}

\usepackage[numbers,sort&compress]{natbib}
\bibpunct[, ]{[}{]}{,}{n}{,}{,}

\theoremstyle{plain}
\newtheorem{theorem}{Theorem}[section]

\newtheorem{corollary}[theorem]{Corollary}

\theoremstyle{definition}

\theoremstyle{remark}

\newcommand{\rP}{\mathrm{P\space }} 
\newcommand{\rE}{\mathrm{E\space }} 

\begin{document}


\title{Theta-positive branching in varying environment}

\author{
\name{S. Sagitov\textsuperscript{a}\thanks{%
CONTACT S. Sagitov. Email: serik@chalmers.se}, A. Lindo \textsuperscript{b}, and Y. Zhumayev \textsuperscript{c}}
\affil{\textsuperscript{a}Chalmers University and University of Gothenburg, 412 96 Gothenburg, Sweden; \textsuperscript{b}University of Glasgow, UK; \textsuperscript{c} L. N. Gumilev Eurasian National University, Kazakhstan}
}

\maketitle

\begin{abstract}
Branching processes in a varying environment encompass a wide range of stochastic demographic models, and their complete understanding in terms of limit behavior poses a formidable research challenge. In this paper, we conduct a thorough investigation of such processes within a continuous-time framework, assuming that the reproduction law of individuals adheres to a specific parametric form for the probability generating function. Our six clear-cut limit theorems support the notion of recognizing five distinct asymptotical regimes for branching in varying environments: supercritical, asymptotically degenerate, critical, strictly subcritical, and loosely subcritical.
\end{abstract}

\begin{keywords}
Continuous time branching process, varying environment, theta branching, limit theorems
\end{keywords}

\section{Introduction}

%
%
%
%
%
%

The subject of this paper is a  time inhomogeneous  Markov branching process $\{Z_{t}\}_{t\ge0}$ with  $Z_{0}=1$. It is a stochastic model for the fluctuating size of a population consisting of individuals that live and reproduce independently of each other, provided that the coexisting individuals are jointly effected by the shared varying environment in the following way:
\begin{description}
\item[ ] - an individual alive at time $t$ dies during the time interval $(t,t+\delta)$ with probability $\lambda_{t}\delta+o(\delta)$ as $\delta\to0$,
\item[ ] - an individual dying at the time $t$ is instantaneously replaced by $k$ offspring with probability $p_{t}(k)$, where $k=0$ or $k\ge2$.
\end{description}
The time-dependent reproduction law of this model is summarized by two functions  
$$ \Lambda_t=\int_0^t \lambda_udu, \quad h_t (s)=p_{t}(0)+p_{t}(2)s^2+p_{t}(3)s^3+\ldots,$$ 
where $h_t (s)$ is the probability generating function for the offspring number and $\Lambda_t$, assumed to be finite for all $t\ge0$, is the cumulative hazard function of the life length of the initial individual. 
In terms of the mean offspring number 
$$a_t=\sum_{n=2}^\infty kp_t(k)=\frac{\partial h_t(s)}{\partial s}|_{s=1},$$also assumed to be finite for all $t\ge0$, the mean population size $\mu_t= \rE(Z_t)$ has the following expression (see Section \ref{4.1})
\begin{equation}\label{MA}
\mu_t=\exp\left\{\int_0^t (a_u-1) d\Lambda_u\right\}.
\end{equation}
Putting $m_t=\rE(Z_t|Z_t>0)$, observe that 
\begin{equation}\label{eext}
\mu_t=m_t\rP(Z_t>0).
\end{equation}

Recall that the extinction probability of the branching process is well defined by
$$q=\lim\rP(Z_t=0).$$
(Here and elsewhere in this paper, the limiting relations are understood to hold  as $t\to\infty$, unless it is clearly stated otherwise.)
In the time homogeneous case, with 
$$\lambda_t\equiv \lambda,\quad a_t\equiv a,\quad \mu_t=e^{(a-1)\lambda t},$$
a Markov branching process  \cite[Ch III]{AN} has one of three possible regimes of reproduction: supercritical when $a>1$, critical when $a=1$, and subcritical when $a<1$. So that in the supercritical case,  $q<1$ and $\mu_t\to \infty$; in the critical case, $q=1$, $\mu_t\equiv 1$, $m_t\to\infty$; and in the subcritical case, $q=1$, $\mu_t\to 0$.  
Compared to the time homogeneous setting, the added feature of varying environment makes the model very flexible and therefore cumbersome to study  
in the most general setting \cite{CM,CH}. In this paper, we distinguish between five classes of the branching processes in varying environment
 \begin{itemize}
\item[(i)] \textit{supercritical} if $q<1$ and $\lim\mu_t=\infty$, 
\item[(ii)]  \textit{asymptotically degenerate} if $q<1$ and $\liminf\mu_t<\infty$,
\item[(iii)] \textit{critical} if $q=1$ and $\lim m_t=\infty$, 
\item[(iv)] \textit{strictly subcritical} if $q=1$ and  $\lim m_t\in[1,\infty)$,
\item[(v)] \textit{loosely subcritical} if $q=1$ and  $\lim m_t$ does not exist.
\end{itemize}

The subject of this paper is a special family of branching processes in varying environment which we call \textit{theta-positive branching process} with the branching parameter $\theta\in(0,1]$ in varying environment $(\{\lambda_t\}, \{a_t\})$. The branching parameter $\theta$  controls the higher moments of the offspring distribution specified by the formula
\begin{equation}\label{pgf}
h_t (s)=1-a_t(1-s)+a_t(1+\theta)^{-1}(1-s)^{1+\theta}.
\end{equation}
It is assumed that the fluctuations of the mean offspring number $a_t$ are restricted to a fixed interval
\begin{equation}\label{keyC}
0\le a_t \le 1+1/\theta.
\end{equation}
This condition guarantees that the probability of zero offspring
$$p_t(0)=h_t(0)=1-(1+\theta)^{-1}\theta a_t$$ 
belongs to the interval $[0,1]$. Observe that given \eqref{keyC}, relation \eqref{MA} implies 
\begin{equation}\label{Ub}
e^{- \Lambda_t}\le \mu_t,\qquad 
\mu_t^\theta\le e^{\Lambda_t}.
\end{equation}

In the important special case of \eqref{pgf} with $\theta=1$, when
\[p_t(0)=1-a_t/2,\quad p_t(2)=a_t/2,\]
the theta-positive branching process turns into the classical birth and death process in varying environment \cite{Ken}.
Such birth-death processes have rich applications in population biology and genetics \cite{St,TN}.
Our study is novel due to the case $0<\theta<1$, where the branching process is featured by the offspring number distribution (see Section \ref{cont})
\begin{align*}
p_t(0)&=1-\theta(1+\theta)^{-1}a_t,\\
p_t(2)&=2^{-1}\theta a_t,\\
p_t(k)&=(k!)^{-1}\theta(1-\theta)(2-\theta)\cdots (k-2-\theta)a_t,\quad k\ge3,
\end{align*}
whose variance is infinite. 
Such theta-positive branching processes might be used in demographic models claiming large variation in the number of offspring. 

The key feature of the theta-positive branching process $Z_t$ is the explicit probability generating function (see Section \ref{4.1})
\begin{equation}\label{thetaF}
\rE(s^{Z_t})=1-(B_{t,\theta}+\mu_t^{-\theta}(1-s)^{-\theta})^{-1/\theta},
\end{equation}
where
\begin{equation}\label{thetaB}
B_{t,\theta}=\theta(1+\theta)^{-1} \int_0^t \mu_u^{-\theta} a_ud\Lambda_u
\end{equation}
is a non-negative term free from the varying $s$, and $\mu_t$ is defined by \eqref{MA}. Notice that with $\theta=1$, the probability generating function  \eqref{thetaF} is a linear-fractional function of $s$. 
In Section \ref{mr} we present the main results of our study based on \eqref{thetaF} and addressing  each of the cases (i)-(v). These results are illustrated in Section \ref{spex} using several worked out special cases and examples. The final Section \ref{cont} contains the proofs.

\subsubsection*{Remarks}
\begin{enumerate}
\item The division into five classes (i)-(v) is a modified version of the classification suggested in \cite{Ke} for the branching processes in varying environment with discrete time.  
In \cite{Ke}, the classes (iv) and (v) are considered as one class called subcritical. 
\item There is a potential for applying the results of this paper in machine learning due to the following recently found link between iterated generating functions and deep neural networks \cite{LT, LPSS}. Consider a fully connected neural network with random weights. Under mild conditions on the activation functions such neural network in the infinite-width limit converges to a Gaussian process \cite{Hanin}. The covariance kernel of this Gaussian process can be calculated in terms of compositions of dual activation functions introduced in~\cite{Daniely}. As it was noted in~\cite{LT}, if $L_{2}$ norm of an activation function with respect to Gaussian measure equals one, then its dual activation is a probability generating function and therefore the corresponding covariance kernel can be expressed using compositions of probability generating functions. 
\item We plan to extend the setting of the current paper using the ideas of \cite{SL,SC} and consider the theta-branching processes in varying environment 
with defective reproduction laws having $h_t(1)<1$. Some inspiration for this future work will come from the recent related paper \cite{KM}.
\item An important direction opened by these results is the study of theta-positive branching processes in random environments; see, for instance, the recent papers \cite{AG, AH} dealing with the discrete-time setting. In light of the latter reference, one might even consider the alternative title \emph{“Power-fractional branching in varying environment”} for this paper.
\end{enumerate}

\section{Main results}\label{mr}

The six theorems of this section deal with a theta-positive branching process in varying environment with parameters $(\theta, \{\lambda_t\}, \{a_t\})$. 
Recall \eqref{MA} and put
\[V_{t,\theta}=\theta(1+\theta)^{-1} \int_0^t \mu_u^{-\theta}d\Lambda_u,\quad V_\theta=\lim V_{t,\theta},\quad \Lambda=\lim \Lambda_t.\]

\begin{theorem}\label{Thm0}
If $V_\theta<\infty$, then $q<1$,  
 \begin{equation}\label{mu1}
\lim \mu_t=\mu,\quad 0<\mu\le\infty,
\end{equation}
and
\begin{equation}\label{qu}
q=1-(V_\theta+(1+\theta)^{-1}+ \theta(1+\theta)^{-1} \mu^{-\theta})^{-1/\theta}.
\end{equation}
If $V_\theta=\infty$, then  $q=1$ and
\begin{align}
\rP(Z_t >0)
&\sim (V_{t,\theta}+\theta(1+\theta)^{-1}\mu_t^{-\theta})^{-1/\theta},  \quad m_t\sim (\mu_t^{\theta}V_{t,\theta}+\theta(1+\theta)^{-1})^{1/\theta}. \label{eex}
\end{align}

\end{theorem}

\begin{theorem}\label{Thm1}
A theta-positive branching process  is supercritical if and only if  $V_\theta<\infty$ and $\Lambda=\infty$. In this case,
 $\lim \mu_t=\infty$,
\begin{equation}\label{qsup}
q=1-(V_\theta+(1+\theta)^{-1})^{-1/\theta},
\end{equation}
and  $\mu_t^{-1}Z_t$ almost surely converges to a random varying $W$ such that
\begin{equation}\label{W}
\rE(e^{-w W})=1-(V_\theta+(1+\theta)^{-1}+w^{-\theta})^{-1/\theta}.
\end{equation}
\end{theorem}

\begin{theorem}\label{Thm2}
A theta-positive branching process  is asymptotically degenerate if and only if $\Lambda<\infty$. In this case, 
 \begin{equation}\label{mu}
\lim \mu_t=\mu,\quad 0<\mu<\infty,
\end{equation}
 and 
$Z_t$ almost surely converges to a random varying $Z_\infty$ such that
\[\rE(s^{Z_\infty}) = 1-(V_\theta +(1+\theta)^{-1}(1-\mu^{-\theta})+\mu^{-\theta} (1-s)^{-\theta})^{-1/\theta}.\]
\end{theorem}

\begin{corollary}
If $\Lambda<\infty$ and $a_t\equiv0$,  then the theta-positive branching process is asymptotically degenerate with $\mu=e^{-\Lambda}$ and $\rE(s^{Z_\infty}) = 1-\mu+\mu s$.
\end{corollary}

\begin{theorem}\label{Thm3}
A theta-positive branching process  is  critical  if and only if 
$V_\theta=\infty$ and 
\begin{equation}\label{gol}
\mu_t^{\theta}V_{t,\theta}\to \infty.
\end{equation}
In this case, 
\begin{align}
 \rP(Z_t >0)&\sim V^{-1/\theta}_{t,\theta}, \qquad 
 m_t\sim \mu_tV^{1/\theta}_{t,\theta},\label{cr2}
  \end{align}
  and
\begin{align}
 \lim\mathcal \rE(e^{-wZ_t/m_t}|Z_t>0)&=1-(1+w^{-\theta})^{-1/\theta},\quad w\ge0. \label{cr3}
  \end{align}

\end{theorem}

\begin{corollary}\label{corr}
If $\Lambda=\infty$ and $0<\liminf \mu_t\le  \limsup  \mu_t<\infty,$  then the theta-positive branching process is critical.
\end{corollary}

\begin{theorem}\label{Thm4}
A theta-positive branching process  is strictly subcritical if and only if $V_\theta=\infty$ and
\begin{equation}\label{go}
 \mu_t^{\theta}V_{t,\theta}\to M_\theta,\quad 0\le M_\theta<\infty.
\end{equation}
In this case, $\mu_t\to0$, 
\begin{align}
& \rP(Z_{t} >0)\sim m^{-1} \mu_t ,  \qquad
m_t\to m,\quad m=(M_\theta+\theta(1+\theta)^{-1} )^{1/\theta},  \label{a3}
\end{align}
and
\begin{equation}\label{a4}
\rE(s^{Z_{t}}|Z_{t}>0)\to1-m(M_\theta- (1+\theta )^{-1}+(1-s)^{-\theta})^{-1/\theta}.
\end{equation}
\end{theorem}
\begin{corollary}
If $\Lambda=\infty$ and 
$ \int_0^\infty a_ud\Lambda_u<\infty,$
then the theta-positive branching process is strictly subcritical.
\end{corollary}

\begin{theorem}\label{Thm5}
A theta-positive branching process  is loosely subcritical if and only if 
$V_\theta=\infty$ and $ \mu_t^{\theta}V_{t,\theta}$ does not have a limit.
In this case, there are several subsequences $t'=\{t_n\}$ leading to different partial limits
\begin{equation}\label{a1}
 \mu_{t'}^{\theta}V_{t',\theta}\to M_\theta,\quad t'\to\infty.
\end{equation}
If \eqref{a1} holds with $M_\theta=\infty$,
then  \eqref{cr2} and \eqref{cr3} are valid with $t=t'$ as  $t'\to\infty$.
On the other hand, if \eqref{a1} holds with $0\le M_\theta<\infty$,
then $\mu_{t'}\to0$, and  \eqref{a3} as well as \eqref{a4} are valid with $t=t'$ as $t'\to\infty$.

\end{theorem}

\noindent\textbf{Remarks}
\begin{enumerate}
\item Theorem \ref{Thm2} describes the well-known asymptotically degenerate case \cite{Lin} when the branching  process in varying environment with a positive probability $1-q$ survives forever as its reproduction process "falls asleep". 
\item Notice that the limiting Laplace transform in  \eqref{cr3}  is the same as in Theorem 7 in \cite{Z} obtained for the critical Markov branching processes in constant environment.
\item For an arbitrary choice of increasing  time points $\{t_n\}$,
the Markov chain $\{Z_{t_n}\}_{n\ge0}$ is a Galton-Watson process in varying environment. Compared to the continuous time setting,
such discrete time branching processes in varying environment
are studied more extensively, see \cite{Ke} and references therein. 
\end{enumerate}

\section{Examples}\label{spex}
Notice that if
$\lambda_t = \lambda(t+1)^\alpha$ and $\alpha<-1$, then $\Lambda<\infty$  implying the asymptotically degenerate case. In the rest of this section we will assume 
\begin{align}
\lambda_t &= \lambda(t+1)^\alpha,\quad 0<\lambda<\infty,\quad -1\le \alpha<\infty. \label{ex1}
\end{align}

\subsection{An example with $a_t\searrow1$}
Assume \eqref{ex1} together with
\begin{align*}
a_t&=1+(1+t)^{-\beta},\quad 0\le \beta<\infty, 
\end{align*}
so that if $\beta\ne 1+\alpha$, then
\[\ln\mu_t
=\lambda(1+\alpha-\beta)^{-1}((1+t)^{1+\alpha-\beta}-1),\]
and if $\beta= 1+\alpha$, then $\mu_t=(1+t)^\lambda$.\\

(a) Suppose $\beta> 1+\alpha$. Then $\Lambda=\infty$ and $\mu_t\to e^{\lambda/(\beta-1-\alpha)}$. This is a critical case according to the Corollary \ref{corr}.\\

(b) Suppose  $\beta= 1+\alpha$. Then $\mu_t$ has a polynomial growth, and 
\[V_\theta=\theta(1+\theta)^{-1}\int_0^\infty \mu_u^{-\theta}d\Lambda_u=\theta(1+\theta)^{-1}\lambda\int_0^\infty (1+u)^{-\theta\lambda+\alpha}du,\]
implying that $V_\theta<\infty$ if and only if $\theta\lambda>1+\alpha$. 
Thus, the case $\lambda^{-1}(1+\alpha)<\theta\le 1$ is supercritical.

If $\theta<\lambda^{-1}(1+\alpha)$, then  $q=1$ and 
\[\mu_t^{\theta} V_{t,\theta}=\theta(1+\theta)^{-1}(1+t)^{\theta\lambda} \int_0^t (1+u)^{-\theta\lambda+\alpha}du\sim \theta(1+\theta)^{-1}(1+\alpha-\theta\lambda)^{-1}(1+t)^{1+\alpha}.\]
This is a critical case since the last relation implies \eqref{gol}.

If $\theta=\lambda^{-1}(1+\alpha)$, then  $q=1$ and we are in the critical case with 
\[\mu_t^{\theta} V_{t,\theta}=\theta(1+\theta)^{-1}(1+t)^{1+\alpha}\ln(1+t).\]

(c) Suppose  $\beta< 1+\alpha$. Then necessarily $\alpha>-1$, $\mu_t\to\infty$, and $V_\theta<\infty$. This is a supercritical case.\\

\noindent\textbf{Remark}. 
According to the item (b),  for a given varying environment $(\{\lambda_t\},\{a_t\})$, the criticality of the theta-positive branching process may depend on the value of the branching parameter $\theta$.

\subsection{An example  with $a_t\nearrow1$}

Assume \eqref{ex1} together with
\begin{align*}
a_t&=1-(1+t)^{-\beta},\quad 0\le \beta<\infty, 
\end{align*}
so that if $\beta\ne 1+\alpha$, then
\[\ln\mu_t=\lambda(1+\alpha-\beta)^{-1}(1-(1+t)^{1+\alpha-\beta}),\]
and if $\beta= 1+\alpha$, then $\mu_t=(1+t)^{-\lambda}$.\\

(a) Suppose  $\beta> 1+\alpha$. Then $\Lambda=\infty$ and $\mu_t\to e^{\lambda/(1+\alpha-\beta)}$. This is a critical case according to the corollary of Theorem \ref{Thm3}.\\

(b) Suppose $\beta= 1+\alpha$. Then  $V_\theta=\infty$ and
\[\mu_t^{\theta} V_{t,\theta}\sim \theta(1+\theta)^{-1} \lambda(1+\alpha+\theta\lambda)^{-1}(1+t)^{1+\alpha}.\]
It follows that the trivial case $\alpha=-1$, $\beta=0$ is strictly subcritical, and the case $\alpha>-1$, $\beta= 1+\alpha$  is critical. \\

(c) Suppose  $\beta< 1+\alpha$. Since 
\[\mu_t=\exp\{\lambda(1+\alpha-\beta)^{-1}(1-(1+t)^{1+\alpha-\beta})\},\]
we find that $V_\theta=\infty$ and $\mu_t^{\theta} V_{t,\theta}$
has a finite limit. This is a strictly subcritical case.

\subsection{An example  with $a_t\searrow0$}

Assume \eqref{ex1} together with
\begin{align*}
a_t&=(1+t)^{-\beta},\quad 0\le \beta<\infty, 
\end{align*}
so that if $\beta\ne 1+\alpha$, then
\[\ln\mu_t=\lambda(1+\alpha-\beta)^{-1}(1-(1+t)^{1+\alpha-\beta})-\lambda(1+\alpha)^{-1}((1+t)^{1+\alpha}-1),\]
and if $\beta= 1+\alpha$, then 
\[\ln\mu_t=\lambda \ln(1+t)-\lambda(1+\alpha)^{-1}((1+t)^{1+\alpha}-1).\]
For this example, $V_\theta=\infty$ and $M_{t,\theta}$ has a finite limit. This is a strictly subcritical case.

\subsection{A loosely subcritical case}\label{vast}

\begin{figure}[h]
\begin{center}
\includegraphics[width=7cm]{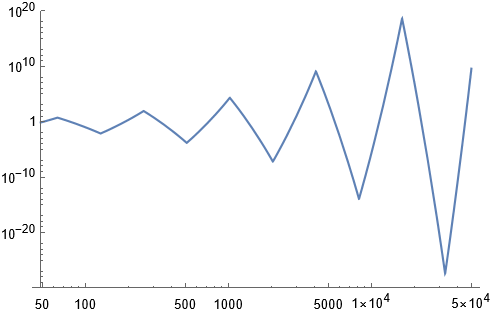}
\caption{On the plot, the horizontal axis gives the time varying $t$ and the vertical axis gives $\mu_t$ for the  example of Section \ref{vast}. 
}
\label{Fig}
\end{center}
\end{figure}

Here we consider a case with vastly alternating environment such that 
\[\liminf\mu_t=0,\quad \limsup\mu_t=\infty.\]
Let $\theta=1$, $\lambda_t=(1+t)^{-1/2}$, and consider a theta-positive branching process with $a_t$ having two alternating values 0 and $2$, so that
\[
a_t-1=\left\{
\begin{array}{rlc}
 1 &  \text{if }  0\le t<2,\\
 -1&   \text{if }   2^{2k-1}\le t <2^{2k} &  \text{ for some }k\ge1,\\
 1 &  \text{if }  2^{2k}\le t <2^{2k+1} &\text{ for some }k\ge1.
\end{array}
\right.
\]
As illustrated by Figure \ref{Fig}, this is a loosely subcritical case with condition \eqref{a1} satisfied for the full range of partial limits $M_\theta\in[0,\infty]$.

\section{Proofs }\label{cont}

We start this section by establishing the properties of   the function \eqref{pgf} stated in the Introduction. To this end, consider
\[h(s)=1-a(1-s)+a(1+\theta)^{-1}(1-s)^{1+\theta}\]
assuming $ \theta\in(0,1]$ and $0\le a\le 1+\theta^{-1}$. Since $h(1)=1$,
$$h(0)=1-a(1+\theta)^{-1}\theta ,\quad h'(0)=0,\quad h''(0)=\theta a, $$
and
\begin{align*}
h^{(k)}(0)&=\theta(1-\theta)(2-\theta)\cdots (k-2-\theta)a,\quad k\ge3,
\end{align*}
are non-negative,
we conclude that
$$h (s)=p(0)+p(2)s^2+p(3)s^3+\ldots,$$ 
is a probability generating function
with 
\begin{align*}
 p(0)&=1-a(1+\theta)^{-1}\theta ,\quad p(1)=0, \quad p(2)=2^{-1}\theta a, \\
 p(k)&=(k!)^{-1}\theta(1-\theta)(2-\theta)\cdots (k-2-\theta)a,\quad k\ge3.
\end{align*}

\subsection{Derivation of \eqref{MA}, \eqref{thetaF}, and \eqref{thetaB}}\label{4.1}
Put for $t\ge\tau$,
\[F_t (\tau,s)=\rE(s^{Z_t}|Z_\tau=1),\quad 0\le s\le1,\]
and notice that $F_t (0,s)=\rE(s^{Z_t})$, see the left hand side of \eqref{thetaF}.
This family of probability generating functions satisfies the following Kolmogorov backward equation, see \cite{Sev},
\begin{equation}\label{kol}
 {\partial F_t (\tau,s)\over\partial \tau}=(F_t (\tau,s)-h_\tau (F_t (\tau,s)))\lambda_\tau,\quad F_t(t,s)=s.
\end{equation}
Setting
\[\mu_t(\tau)={\partial F_t (\tau,s)\over\partial s}|_{s=1},\]
we find from \eqref{kol}
\[ {\partial \mu_t(\tau)\over\partial \tau}=\mu_t(\tau)(1-a_\tau)\lambda_\tau,\quad \mu_t(t)=1\]
implying 
$$\mu_t(\tau)=\exp\Big\{\int_\tau^t (a_u-1) d\Lambda_u\Big\}.$$ 
Setting $\tau=0$ in the last expression, we arrive at \eqref{MA}.

With the generating function $h_t(s) $ having the special form \eqref{pgf}, the equation \eqref{kol} yields a Bernoulli differential equation for $x_\tau =1-F_t (\tau,s)$,
 \[x_\tau '=(1-a_\tau)\lambda_\tau x_\tau +(1+\theta)^{-1}a_\tau\lambda_\tau x_\tau ^{1+\theta},\quad x_t=1-s.\]
In terms of $y_\tau =x_\tau ^{-\theta}$ it leads through $y_\tau '=-\theta x_\tau 'x_\tau^{-1-\theta} $  to a linear differential equation
 \[y_\tau '=\theta (a_\tau-1) \lambda_\tau y_\tau - \theta (1+\theta)^{-1}a_\tau\lambda_\tau,\quad y_t=(1-s)^{-\theta},\]
 which has an explicit solution
\[y_\tau =\mu_\tau^{\theta}(\mu_t^{-\theta }(1-s)^{-\theta}+B_{t,\theta}-B_{\tau,\theta}), \]
where $B_{t,\theta}$ is given in \eqref{thetaB}. 
Furthermore,  in view of $F_t (\tau,s)=1-y_\tau^{-1/\theta}$, we conclude that
\begin{align}
 F_t (\tau,s)&=1-\mu_\tau^{-1}(\mu_t^{-\theta}(1-s)^{-\theta}+B_{t,\theta}-B_{\tau,\theta})^{-1/\theta}. \label{AS}
\end{align}
Setting $\tau=0$ in the last relation we arrive at
\eqref{thetaF} with  \eqref{thetaB}. 

\subsection{Proof of Theorem \ref{Thm0}}\label{4.2}

Observe that in view of \eqref{MA} and \eqref{thetaB}, the derivative over $t$
\[(B_{t,\theta}+\mu_t^{-\theta})' 
=\theta \mu_t^{-\theta}\lambda_t(1-\theta(1+\theta)^{-1}a_t)\]
is non-negative due to \eqref{keyC}.
By integration,
\begin{align*}
B_{t,\theta}+\mu_t^{-\theta}-1=(1+\theta) V_{t,\theta}-\theta B_{t,\theta},
\end{align*}
entailing 
\begin{align}\label{NY}
B_{t,\theta}=V_{t,\theta}+(1+\theta)^{-1}(1-\mu_t^{-\theta}).
\end{align}
Observe also that setting $s=0$ in \eqref{thetaF} gives
\begin{align*}
\rP(Z_t >0)&=(B_{t,\theta}+\mu_t^{-\theta})^{-1/\theta}. 
\end{align*}
This together with relations  \eqref{thetaF} and \eqref{NY} yield
\begin{align}
\rE(s^{Z_t})&=1-(V_{t,\theta}+(1+\theta)^{-1}(1-\mu_t^{-\theta})+\mu_t^{-\theta}(1-s)^{-\theta})^{-1/\theta}, \label{NYe}\\
\rP(Z_t >0)&=(V_{t,\theta}+(1+\theta)^{-1}+\theta(1+\theta)^{-1}\mu_t^{-\theta})^{-1/\theta}. \label{NYs}
\end{align}

Turning to the statement of  Theorem \ref{Thm0}, assume first that $V_\theta<\infty$. By \eqref{thetaB} and \eqref{keyC}, the limit $B_{\theta}=\lim B_{t,\theta}$ always exists and  $B_{\theta}\le V_\theta$. 
Therefore, relation \eqref{NY} implies the existence of $\lim\mu_t=\mu$ for some $\mu\in(0,\infty]$. Combing this with \eqref{NYs} we conclude that $q$ satisfies \eqref{qu}, so that $q<1$.

On the other hand, given $V_\theta=\infty$, relations \eqref{NYs} and \eqref{eext} imply  $q=1$ together with  \eqref{eex}.

\subsection{Proof of Theorem \ref{Thm1}}
In the supercritical case,  when $q<1$ and $\mu_t\to\infty$, Theorem \ref{Thm0} gives $V_\theta<\infty$, relation \eqref{Ub} implies $\Lambda=\infty$, and relation \eqref{qsup} follows from \eqref{qu}. On the other hand, if $V_\theta<\infty$ and $\Lambda=\infty$, then by Theorem \ref{Thm0}, we have $q<1$ and \eqref{mu1}. From \eqref{mu1}, $\Lambda=\infty$, and $V_\theta<\infty$, we derive $\mu=\infty$.

 As a non-negative martingale, $W_t=\mu_t^{-1}Z_t$  almost surely converges to a limit $W$, and it remains to prove relation \eqref{W}. The  Laplace transform of $W_t$ computed using \eqref{thetaF} gives
\[\rE(e^{-wW_t})=1-(B_{t,\theta}+(\mu_t(1-e^{-w /\mu_t}))^{-\theta})^{-1/\theta}\to 1-(B_{\theta}+w^{-\theta})^{-1/\theta},\]
which together with  \eqref{NY} yields \eqref{W}. 

\subsection{Proof of Theorem \ref{Thm2}  }
Put
$$A_t=\int_0^ta_ud\Lambda_u,\quad A=\lim A_t.$$
 If $\Lambda<\infty$, then $A<\infty$ due to the condition \eqref{keyC}, and \eqref{mu} holds with $\mu=e^{A-\Lambda}$.
 This entails $V_\theta<\infty$, so that according to Theorem \ref{Thm1}, we have  $q<1$ and we are in the asymptotically degenerate case. 
On the other hand, if we are in the asymptotically degenerate case, then by Theorem \ref{Thm0}, we have $V_\theta<\infty$. According to Theorem \ref{Thm1}, the relation $\Lambda=\infty$ would imply the supercritical case, thus we must have $\Lambda<\infty$. 
We conclude that $\Lambda<\infty$ is a necessary and sufficient condition for the asymptotically degenerate case. 

Assume  $\Lambda<\infty$ and put 
\[P(\tau,t)=\rP(Z_t=1|Z_\tau=1),\quad 0\le\tau\le t.\]To prove the stated almost sure convergence it suffices to verify  the following Lindvall's condition \cite{Lin}:
 \begin{equation}\label{church}
 \sum_{n\ge1}(1-P(t_n,t_{n+1}))<\infty
\end{equation}
for an arbitrary  sequence $\{t_n\}$ monotonely increasing to infinity. Taking the derivative over $s$ in \eqref{AS}
\[ {\partial F_t (\tau,s)\over\partial s}=\mu_\tau^{-1}\mu_t^{-\theta}(1-s)^{-\theta-1}(\mu_t^{-\theta}(1-s)^{-\theta}+B_t(\theta)-B_\tau(\theta))^{-1/\theta-1} \]
and setting here $s=0$  we get 
\[P(\tau,t)=\mu_\tau^{-1}\mu_t(1+\mu_t^{\theta}(B_t(\theta)-B_\tau(\theta)))^{-1/\theta-1}.\]
We prove \eqref{church} by using the upper bounds
\begin{align*}
1-P(\tau,t)&\le (1/\theta+1)\mu_\tau^{-1}\mu_t^{1+\theta}(B_t(\theta)-B_\tau(\theta))+1-\mu_\tau^{-1}\mu_t\\
1-\mu_\tau^{-1}\mu_t&=1-e^{A_t-A_\tau}e^{\Lambda_\tau-\Lambda_t}\le1-e^{\Lambda_\tau-\Lambda_t}\le \Lambda_t-\Lambda_\tau.
\end{align*}
These bounds together with \eqref{mu} imply the existence of a positive constant $c$, such that
\begin{align*}
 \sum_{n\ge1}(1-P(t_n,t_{n+1}))&\le cB_\theta+\Lambda<\infty.
\end{align*}

\subsection{Proof of Theorem \ref{Thm3} }

The stated criticality conditions $V_\theta=\infty$ and \eqref{gol} as well as \eqref{cr2}  immediately follow from Theorem \ref{Thm0}. Relation \eqref{cr3}, is verified using
  \begin{equation}\label{va}
  \rE(s^{Z_t}|Z_t>0)=\frac{\rE(s^{Z_t})-\rP(Z_t =0)}{\rP(Z_t >0)}=1-\frac{1-\rE(s^{Z_t})}{\rP(Z_t >0)}.
\end{equation}
Applying \eqref{NYe} and \eqref{NYs} we obtain
 \[\frac{1-\rE(e^{-wZ_t/m_t})}{\rP(Z_t >0)}=\frac{(1+(1+\theta)^{-1}V_{t,\theta}^{-1}(1-\mu_t^{-\theta})+V_{t,\theta}^{-1}\mu_t^{-\theta}(1-e^{-w/m_t})^{-\theta})^{-1/\theta}}
 {(1+(1+\theta)^{-1}V_{t,\theta}^{-1}+\theta(1+\theta)^{-1}V_{t,\theta}^{-1}\mu_t^{-\theta})^{-1/\theta}}.   \]
As $t\to\infty$, this together with \eqref{gol}, \eqref{cr2}, and  \eqref{va} yield \eqref{cr3}.

\subsection{Proof of Theorem \ref{Thm4} and Theorem \ref{Thm5} }
The stated strict sub-criticality conditions  follow from Theorems \ref{Thm0} and \ref{Thm3}.  Relations $V_\theta=\infty$ and \eqref{go} imply $\mu_t\to0$. The statement  \eqref{a3} follows from \eqref{eex}. The convergence  \eqref{a4} is easily derived from  \eqref{va}. This finishes the proof of Theorem \ref{Thm4}.

Theorem \ref{Thm5} follows from Theorems \ref{Thm3} and \ref{Thm4}. 

\section*{Acknowledgements}
This project was partially supported by the travel grant MF2022-0001 from the GS Magnusons Fond.
The corresponding author is grateful to an anonymous reviewer for a careful reading and for helping to correct some of the stated formulas.

\end{document}